\documentclass{amsproc}%
\usepackage{graphicx}
\usepackage{amscd}
\usepackage{amsmath}
\usepackage{amsfonts}
\usepackage{amssymb}%
\theoremstyle{plain}

\numberwithin{equation}{section}

\begin{document}
\title[The Skolem-Noether Theorem ]{A Constructive Elementary Proof of the\\Skolem-Noether Theorem for Matrix Algebras}
\author{Jen\H{o} Szigeti}
\address{\noindent\noindent Institute of Mathematics, University of Miskolc, Miskolc,
Hungary 3515}
\email{jeno.szigeti@uni-miskolc.hu}
\author{Leon van Wyk}
\address{Department of Mathematical Sciences, Stellenbosch University, P/Bag X1,
Matieland 7602, Stellenbosch, South Africa }
\email{LvW@sun.ac.za}
\thanks{The first author was partially supported by the National Research, Development
and Innovation Office of Hungary (NKFIH) K119934.}
\thanks{The second author was supported by the National Research Foundation of South
Africa under Grant No.~UID 72375. Any opinion, findings and conclusions or
recommendations expressed in this material are those of the authors and
therefore the National Research Foundation does not accept any liability in
regard thereto.}

\begin{abstract}
We give a constructive elementary proof for the fact that any $K$-automorphism
$\varphi$ of the full $n\times n$ matrix algebra over a field $K$ is
conjugation by some invertible $n\times n$ matrix $A$ over $K$.

\end{abstract}
\maketitle

\noindent The theorem stating that any two embeddings of an extension of a
field $K$ into a finite-dimensional central simple algebra over $K$ are
conjugate was first published by Skolem in 1927 (see [Sk]). This theorem was
rediscovered by Noether in 1933 (see [N]). The Skolem-Noether theorem is now
considered as a fundamental result in the theory of central simple algebras,
linear groups and representation theory. See, for example, [GSz] and [R].

\noindent A short constructive proof of the following form of the
Skolem-Noether theorem for matrix algebras can be found in [Se].

\bigskip

\noindent\textbf{Theorem.}\textit{ Any }$K$\textit{-automorphism }%
$\varphi:\mathrm{M}_{n\times n}(K)\longrightarrow\mathrm{M}_{n\times n}%
(K)$\textit{ of the full }$n\times n$\textit{ matrix algebra over a field }%
$K$\textit{ is conjugation by some invertible matrix }$A\in\mathrm{M}_{n\times
n}(K)$\textit{, i.e., }$\varphi(X)=AXA^{-1}$\textit{ for all }$X\in
\mathrm{M}_{n\times n}(K)$\textit{.}

\bigskip

\noindent Although there are many alternative proofs of the Skolem-Noether
theorem in the literature, they are generally not constructive. The value of
the present note is that it provides an elementary proof of the theorem above,
which is entirely constructive and gives the matrix $A$ explicitly, using the
$\varphi$-images of only two matrices (irrespective of the value of $n$), a
nonzero vector in a certain kernel, and matrix multiplication. The proof in
[Se] is different from ours and requires the images of $n$ matrices.

\bigskip

\noindent Our procedure comprises the following standard operations:

\bigskip

\noindent(1) For $1\leq i,j\leq n$, let $E_{i,j}$ denote the matrix in
$\mathrm{M}_{n\times n}(K)$\ having $1$ in the $(i,j)$ position and zeros in
all other positions, and let $S=E_{1,2}+E_{2,3}+\cdots+E_{n-1,n}$. In other
words,%
\[
E_{n,1}=\left[
\begin{array}
[c]{ccccc}%
0 & 0 & \cdots & \cdots & 0\\
\vdots & \vdots &  &  & \vdots\\
\vdots & \vdots &  &  & \vdots\\
0 & 0 & \cdots & \cdots & 0\\
1 & 0 & \cdots & \cdots & 0
\end{array}
\right]  \text{ and }S=\left[
\begin{array}
[c]{ccccc}%
0 & 1 & 0 & \cdots & 0\\
\vdots & \ddots & \ddots & \ddots & \vdots\\
\vdots &  & \ddots & \ddots & \vdots\\
\vdots &  &  & \ddots & 1\\
0 & \cdots & \cdots & \cdots & 0
\end{array}
\right]  .
\]
\noindent(2) Let $H=\varphi(E_{n,1})$ and $G=\varphi(S)$, and let $\mathbf{a}$
be a nonzero column vector in the kernel of $I_{n}-G^{n-1}H=I_{n}%
-\varphi(E_{1,1})$, where $I_{n}$ is the $n\times n$ identity matrix.

\noindent(3) Take the matrix%
\[
A=\left[  G^{n-1}H\mathbf{a}\mid G^{n-2}H\mathbf{a}\mid\cdots\mid
GH\mathbf{a}\mid H\mathbf{a}\right]
\]
with column vectors $G^{n-i}H\mathbf{a}$, $1\leq i\leq n$.

\bigskip

\noindent\textbf{Proof.} Observe that%
\[
I_{n}-\varphi(E_{1,1})=\varphi(I_{n}-E_{1,1})=\varphi(E_{2,2}+E_{3,3}%
+\cdots+E_{n,n})
\]
is not invertible, otherwise $E_{2,2}+E_{3,3}+\cdots+E_{n,n}$ would also be
invertible. Since $\det(I_{n}-\varphi(E_{1,1}))=0$, there exists an
$\mathbf{a}\in\mathrm{M}_{n\times1}(K)$, $\mathbf{a}\neq\mathbf{0}$, such that
$(I_{n}-\varphi(E_{1,1}))\mathbf{a}=\mathbf{0}$. Clearly, $S^{n-1}=E_{1,n}$
and $S^{n-1}E_{n,1}=E_{1,1}$ give that $G^{n-1}H=\varphi(E_{1,1})$ and
$G^{n-1}H\mathbf{a}=\mathbf{a}$. Notice that $S^{i}E_{j,1}=E_{j-i,1}$ for
$1\leq i<j\leq n$, and so $S^{n}=E_{n,1}S^{n-2}E_{n,1}=\cdots=E_{n,1}%
SE_{n,1}=E_{n,1}^{2}=0$ implies $G^{n}=HG^{n-2}H=\cdots=HGH=H^{2}=0$.

\noindent We claim that $AE_{n,1}=HA$ and $AS=GA$, where%
\[
A=\left[  G^{n-1}H\mathbf{a}\mid G^{n-2}H\mathbf{a}\mid\cdots\mid
GH\mathbf{a}\mid H\mathbf{a}\right]
\]
is the matrix in $\mathrm{M}_{n\times n}(K)$ with column vectors
$G^{n-i}H\mathbf{a}$, $1\leq i\leq n$. Since%
\[
AE_{n,1}=\left[  H\mathbf{a}\mid\mathbf{0}\mid\cdots\mid\mathbf{0}%
\mid\mathbf{0}\right]
\]
and%
\[
HA=\left[  HG^{n-1}H\mathbf{a}\mid HG^{n-2}H\mathbf{a}\mid\cdots\mid
HGH\mathbf{a}\mid H^{2}\mathbf{a}\right]  ,
\]
the fact that $G^{n-1}H\mathbf{a}=\mathbf{a}$\ and $HG^{n-2}H=\cdots
=HGH=H^{2}=0$ gives that $AE_{n,1}=HA$. Since%
\[
AS=AE_{1,2}+AE_{2,3}+\cdots+AE_{n-1,n}=\left[  \mathbf{0}\mid G^{n-1}%
H\mathbf{a}\mid\cdots\mid G^{2}H\mathbf{a}\mid GH\mathbf{a}\right]
\]
and%
\[
GA=\left[  G^{n}H\mathbf{a}\mid G^{n-1}H\mathbf{a}\mid\cdots\mid
G^{2}H\mathbf{a}\mid GH\mathbf{a}\right]  ,
\]
the fact that $G^{n}=0$ gives that $AS=GA$. To prove that $A$ is invertible,
take a zero linear combination of the column vectors of $A$:%
\[
\lambda_{n-1}G^{n-1}H\mathbf{a}+\lambda_{n-2}G^{n-2}H\mathbf{a}+\cdots
+\lambda_{1}GH\mathbf{a}+\lambda_{0}H\mathbf{a}=\mathbf{0}.
\]
In view of \ $G^{n-1}H\mathbf{a}=\mathbf{a}\neq\mathbf{0}$\ and $G^{n}=0$,
left multiplication by $G^{n-1}$ gives that $\lambda_{0}\mathbf{a}=\lambda
_{0}G^{n-1}H\mathbf{a}=\mathbf{0}$, whence $\lambda_{0}=0$ follows. Then left
multiplication by $G^{n-2}$ gives $\lambda_{1}\mathbf{a}=\lambda_{1}%
G^{n-1}H\mathbf{a}=\mathbf{0}$, whence $\lambda_{1}=0$ follows. Repeating the
left multiplications, we obtain that $\lambda_{0}=\lambda_{1}=\cdots
=\lambda_{n-1}=0$. Thus the columns of $A$\ are linearly independent, whence
the invertibility of $A$ follows.

\noindent Now $AE_{n,1}=HA$ and $AS=GA$ imply that $AE_{n,1}A^{-1}%
=\varphi(E_{n,1})$ and $ASA^{-1}=\varphi(S)$. Since $E_{i,j}=S^{n-i}%
E_{n,1}S^{j-1}$ for all $1\leq i,j\leq n$, the matrices $E_{n,1}$\ and
$S$\ generate $\mathrm{M}_{n\times n}(K)$ as a $K$-algebra, and so
$AXA^{-1}=\varphi(X)$ for all $X\in\mathrm{M}_{n\times n}(K)$. $\square$

\bigskip

\noindent Although this proof does not directly translate to the more general
case, the theorem immediately implies the statement for central simple
algebras by invoking Hilbert's Theorem 90. That is how Skolem-Noether is
proved in [GSz], however [GSz] gives a nonconstructive proof of the theorem
about matrices.

\noindent The advantage of our proof is that it may help to characterize the
automorphisms of certain two-generator algebras, where the generators satisfy
some of the relations satisfied by our two simple generating matrices. Put
more precisely, our methods may help to determine the $K$-automorphisms of the
quotient of the two-generated free associative (polynomial) algebra
$K<x,y>/N$, where the ideal $N$\ is generated by the monomials $x^{n}%
,yx^{n-2}y,\ldots,yxy,y^{2}$ (perhaps under some additional conditions).

\bigskip

\noindent REFERENCES

\bigskip

\noindent\lbrack GSz] Gille, Ph.; Szamuely, T. \textit{Central simple algebras
and Galois cohomology}, Cambridge Studies in Advanced Mathematics 101,
Cambridge University Press, 2006.

\noindent\lbrack N] Noether, E. \textit{Nichtkommutative Algebra}, Math. Z. 37
(1933), 514-541.

\noindent\lbrack R] Rowen,\ L. H.\textit{ Ring Theory, Vol. II,} Academic
Press, New York, 1988.

\noindent\lbrack Se] Semrl, P. \textit{Maps on matrix spaces}, Linear Algebra
Appl. 413 (2006), 364-393.

\noindent\lbrack Sk] Skolem, T. \textit{Zur Theorie der assoziativen
Zahlensysteme}, Skrifter Oslo 12 (1927), 50.

\end{document}